\documentclass[12pt]{article}

\usepackage{amsmath}
\allowdisplaybreaks[4]

\usepackage[all]{xy}
\usepackage{amssymb}
\usepackage{amsthm}
\usepackage{hyperref}
\hypersetup{colorlinks=true,linkcolor=blue,citecolor=red}
\usepackage{amsmath}
\usepackage{amscd}
\usepackage{verbatim}
\usepackage{eurosym}
\usepackage{float}
\usepackage{color}
\usepackage{dcolumn}
\usepackage[mathscr]{eucal}
\usepackage[all]{xy}
\usepackage{hyperref}
\usepackage{mathrsfs}
\usepackage{amsmath}
\usepackage{amssymb}
\usepackage{amsfonts,ifpdf}
\usepackage{graphicx}
\usepackage{times}
\usepackage{float}
\usepackage{epstopdf}
\usepackage{cite}
\usepackage{youngtab}
\usepackage{ytableau}
\ytableausetup
{mathmode, boxsize=0.9em}

\newtheorem{theo}{Theorem}[section]

\newtheorem{lem}[theo]{Lemma}

\newtheorem{pf}{proof}[section]
\theoremstyle{remark}

\makeatletter \@addtoreset{equation}{section} \makeatother
\makeindex 
\def\qed{\hfill \rule{4pt}{7pt}}

\begin{document}

\begin{center}
{\Large\bf Semi-invariants of Binary Forms and \\[6pt]
Sylvester's Theorem}\\ [7pt]
\end{center}

\vskip 3mm

\begin{center}
William Y.C. Chen$^{1}$ and Ivy D.D. Jia$^{2}$\\[8pt]
$^{1}$Center for Applied Mathematics\\
Tianjin University\\
Tianjin 300072, P. R. China\\[12pt]

$^{2}$School of Science
\\
Tianjin University of Commerce\\
Tianjin 300134, P. R. China\\[15pt]

Emails: $^{1}$chenyc@tju.edu.cn, $^{2}$ivyjia@tjcu.edu.cn
\\[15pt]
\textit{Dedicated to Doron Zeilberger on the Occasion of His Seventieth Birthday} \\[15pt]

\end{center}

\vskip 3mm

\begin{abstract}

We obtain a combinatorial
formula related to the shear transformation for semi-invariants
of binary forms, which implies the classical characterization of   semi-invariants
 in terms of a differential operator.
 Then, we present a combinatorial proof of an identity of Hilbert,
 which leads to a relation of Cayley on semi-invariants. This identity plays a crucial role
in the original proof of Sylvester's theorem on semi-invariants in connection with
the Gaussian coefficients. Moreover, we show that the additivity
lemma of Pak and Panova which yields the strict unimodality of the
Gaussian coefficients for $n,k \geq 8$ can be deduced from the ring property of semi-invariants.

\vskip 6pt

\noindent
{\bf Mathematics Subject Classification:} 05A17, 05E10, 13A50
\\ [7pt]
{\bf Keywords:} Hilbert's identity,  Sylvester's theorem, binary forms, semi-diagrams, semi-invariants, Gaussian coefficients, strict unimodality
\end{abstract}

\section{Introduction}

This work is a continuation of the exploration of
the Gaussian coefficients or the $q$-binomial
coefficients by means of semi-invariants of
binary forms recently carried out in \cite{CJ20}.
We will be mainly concerned with the combinatorial perspectives related to
 Sylvester's proof of the unimodality conjecture of Cayley.
 A key identity
 used by Sylvester is a relation due to Cayley, which turns out to be a
 consequence of an identity of Hilbert.
 We shall give a combinatorial interpretation of the identity of Hilbert.
 Moreover, we show that the additivity lemma of Pak and Panova leading to
 the strict unimodality of the Gaussian coefficients for $n,k\geq 8$ can
 be deduced from the ring property of semi-invariants.

Let
$p(k,n,m)$ denote the number of partitions of $m$
 contained in a $k\times n$ rectangle, then
the Gaussian coefficients can be expressed as
\begin{equation}
{n+k \brack k } =\sum_{m=0}^{nk} p(k,n,m) q^m,
\end{equation}
 see \cite{Andrews98,Stanley12}.

The Gaussian coefficients are symmetric in $q$.
Cayley \cite{Cayley56}  conjectured in 1856 that the Gaussian coefficients are unimodal, and
 it was proved by Sylvester \cite{Sylvester78b} in 1878
 resorting to semi-invariants of binary forms. There has been an extensive literature on this subject,
see, for example, \cite{O'Hara90,PP14,Proctor82,Stanley80,White80,Zeilberger89}.
It is worth mentioning that  O'Hara \cite{O'Hara90} found a constructive proof.
 Zeilberger \cite{Zeilberger89} discovered an identity, known as the KOH theorem, which justifies the unimodality.

As the first step to bring Sylvester's theorem to a combinatorial ground,
let us take a look at the classical characterization of semi-invariants in
terms of the differential operator $D$.
A semi-invariant
of a binary $n$-form is a polynomial $I(a_0, a_1, \ldots, a_n)$ with rational coefficients
such that
\begin{equation}\label{II-1}
I(a_0, a_1, \ldots, a_n) = I(a'_0, a'_1, \ldots, a'_n),
\end{equation}
where the $a_i'$ are determined by the shear transformation with respect to a variable $z$, that is, for $0\leq i\leq n$,
\begin{equation}
a_i'=a_i+{i\choose1}a_{i-1}z+{i\choose2} a_{i-2}z^2+\cdots+a_0 z^i.
\end{equation}

Let
\begin{equation} \label{D}
 {D} = a_0\frac{\partial}{\partial a_1}+2a_1\frac{\partial}{\partial a_2}
+3a_2\frac{\partial}{\partial a_3}+\cdots+na_{n-1}\frac{\partial}{\partial a_n}.
\end{equation}
Semi-invariants can be characterized in terms of the   differential operator $D$, see Cayley \cite{Cayley56},
Sylvester \cite{Sylvester78b}, or Hilbert \cite{Hilbert93}. More precisely,
a polynomial $I(a_0,a_1,\ldots,a_n)$
is a semi-invariant of a binary $n$-form if and only if   $ {D}(I)=0$.

Note that if $I$ and $J$
are two semi-invariants of a binary $n$-form, then so are $I+J$ and $IJ$.
Based on a  combinatorial interpretation of the operator $D$, we show that
 $I(a_0',a'_1,\ldots,a'_n)$ can be expressed in terms
of the polynomial $I(a_0, a_1, \ldots, a_n)$ and the operator $D$, in the spirit of the Taylor expansion. This formula immediately leads to the
characterization of semi-invariants in terms of the operator $D$.

The second objective of this paper is to
present a combinatorial proof of an
identity  of Hilbert involving the operators
$D$ and $\Delta$. The operator $\Delta$ is defined by
\begin{equation}\label{Delta}
\Delta=na_1\frac{\partial}{\partial a_0}+(n-1)a_2\frac{\partial}{\partial a_1}+\cdots+a_n\frac{\partial}{\partial a_{n-1}}.
\end{equation}

The identity of Hilbert \cite{Hilbert93} reads as follows:
For $n,k\geq 0$, and $0\leq m\leq nk$, let $\lambda$ be a partition of $m$ contained in a $k\times n$ rectangle, and let $c=nk-2m$. Then for $i\geq1$,
\begin{equation}\label{kf}
D\Delta^i(a_\lambda)-\Delta^i D(a_\lambda)=i(c-i+1)\Delta^{i-1}(a_\lambda).
\end{equation}

As an application of Sylvester's theorem, we show that the
additivity lemma of Pak and Panova \cite{PP13a,PP13b} can be deduced from the
ring property of semi-invariants. The additivity lemma for the Gaussian
 coefficients was established via a connection with
the Kronecker coefficients
in the representation theory of the symmetric group as well as
the semigroup property of the Kronecker coefficients due to Christandl,  Harrow  and  Mitchison \cite{CHM07}.

\begin{lem}  \label{add}
Assume that $k_1, k_2, n \geq 2$,
    at least one of $k_1$, $k_2$ and $n$ is greater than two and at least
    one of $k_1$, $k_2$ and $n$ is even. If
  the strict unimodality holds for   ${n+k_1\brack n}$ and ${n+k_2\brack n}$,
  then it holds for ${n+k_1+k_2\brack n}$.
\end{lem}

To conclude the introduction, we recall that the strict unimodality proved by Pak and Panova says that for $n,k\geq 8$ and $2\leq m\leq nk/2$,
\begin{equation}
p(k,n,m)> p(k,n,m-1).
\end{equation}

\section{The Operators $D$ and $\Delta$}

A binary form of degree $n$, or a binary $n$-form, is a homogeneous polynomial
in $x$ and $y$,
\begin{equation} \label{F-1}
f(x,y)= a_0x^n+ \binom n1 a_{1} x^{n-1}y+{n\choose 2} a_{2} x^{n-2}y^2+\cdots+a_ny^n,
\end{equation}
where the coefficients $a_0, a_1, \ldots, a_n$ are regarded as
variables.
Consider the shear transformation:
 $x= x' + z y'$ and $y=y'$, where
 $z$ is treated as a variable.
Suppose that under this transformation, the binary form $f(x,y)$ becomes
\begin{align}\label{F-2}
f'(x',y')&=a_0'{x'}^n+{n\choose1}a_1'{x'}^{n-1}y'+\cdots+a_n'{y'}^n.
\end{align}
It is easily checked that for $0\leq i\leq n$,
\begin{equation}\label{ai}
a_i'=a_i+{i\choose1}a_{i-1}z+{i\choose2} a_{i-2}z^2+\cdots+a_0 z^i.
\end{equation}
We say that
a polynomial $I(a_0, a_1, \ldots, a_n)$ with rational coefficients
is a semi-invariant
of the binary form $f(x,y)$ if
\begin{equation}\label{II}
I(a_0, a_1, \ldots, a_n) = I(a'_0, a'_1, \ldots, a'_n),
\end{equation}
 see, for example, \cite{Grosshans03}.

To determine whether a polynomial $I(a_0, a_1, \ldots, a_n)$ satisfies the above condition
\eqref{II}, we are led to an expansion of $I(a'_0, a'_1, \ldots, a_n')$ as a polynomial
of $z$. As expected, the operator $D$ comes to the scene.
Recall that $D$ is given by
\begin{equation} \label{D2}
 {D} = a_0\frac{\partial}{\partial a_1}+2a_1\frac{\partial}{\partial a_2}
+3a_2\frac{\partial}{\partial a_3}+\cdots+na_{n-1}\frac{\partial}{\partial a_n}.
\end{equation}

\begin{theo}\label{SI}
For $n\geq 0$, and for
any polynomial $I(a_0, a_1, \ldots, a_n)$ over the rational numbers, we have
\begin{equation}\label{Sz}
I(a'_0,a'_1,\ldots,a'_n)=\sum_{i\geq0}D^iI(a_0, a_1, \ldots, a_n)\frac{z^i}{i!}.
\end{equation}
\end{theo}

While we do not intend
to claim that the above formula \eqref{Sz} is new due to the lack of accessible literature, at least it is worth noting that such a formulation makes
the characterization of semi-invariants
transparent in the sense that if $D(I)$ vanishes, then so does $D^i(I)$ for any
$i\geq 2$. As will be seen, the idea behind \eqref{Sz}
serves as an embarkation point to
a combinatorial understanding of the identity \eqref{kf}
of Hilbert.

Here is an illustration of Theorem \ref{SI}.
For a monomial $a^\nu = a_0^{\nu_0} a_1^{\nu_1}\cdots a_n^{\nu_n}$, we define
its degree by
\[ \nu_0 + \nu_1+\cdots +\nu_n\]
 and its weight by
\[ \nu_1+2 \nu_2+\cdots+n \nu_n.\]
It is clear that when the operator $D$ is applied to a monomial
$a^\nu$, it preserves the degree and decreases the weight by one.
For example, let \[ f(x,y)=a_0x^3+3a_1x^2y+3a_2xy^2+a_3y^3,\]
 and let
\begin{equation}
I(a_0,a_1,a_2,a_3)=
c_1a_0^2a_3+c_2a_0a_1a_2+c_3a_1^3.
\end{equation}
Upon the substitution \eqref{ai}, we get
\[
\begin{split}
 I(a'_0,a'_1,a'_2,a'_3)&=c_1a_0^2a_3+c_2a_0a_1a_2+c_3a_1^3\\[5pt]
 &\quad +((3c_1+c_2)a_0^2a_2+(2c_2+3c_3)a_0a_1^2)z\\[5pt]
 &\quad +3(c_1+c_2+c_3)a_0^2a_1z^2\\[5pt]
 &\quad +(c_1+c_2+c_3)a_0^3z^3.
\end{split}
\]
On the other hand,
\[
\begin{split}
D(I)&=(3c_1+c_2)a_0^2a_2+(2c_2+3c_3)a_0a_1^2,
\\[5pt]
D^2(I)&=6(c_1+c_2+c_3)a_0^2a_1,
\\[5pt]
D^3(I)&=6(c_1+c_2+c_3)a_0^3.
\end{split}
\]
We see that \eqref{Sz} holds in this case.

To lay out a combinatorial setting for Theorem \ref{SI},
we adopt the common notation of a partition $\lambda$ of $m$ contained
in a $k \times n$ rectangle, that is, $\lambda=(\lambda_1, \lambda_2, \ldots, \lambda_k)$,
where $n\geq \lambda_1 \geq \lambda_2 \geq \cdots \geq \lambda_k \geq 0$ and
$\lambda_1+\lambda_2+\cdots + \lambda_k=m$. Accordingly, we use $a_\lambda$ to
denote the monomial $a_{\lambda_1}a_{\lambda_2}\cdots a_{\lambda_k}$.

\begin{figure}
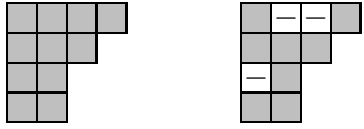

\centering
\ytableaushort{\ \ \ \ ,\ \ \ ,\ \ ,\ \ }
* [*(lightgray)]{4,3,2,2}
\hspace{40pt}
\ytableaushort{\ {*(white)-} {*(white)-} \ ,\ \ \ ,{*(white)-}\ ,\ \ }* [*(lightgray)]{4,3,2,2}
\vspace{0cm}
\caption{The Young diagram and a semi-diagram of shape $(4,3,2,2)$.}
\label{4322}
\end{figure}

Recall that the Young diagram or the Ferrers diagram of
a partition $\lambda$ is a collection of shaded boxes or cells arranged in left-justified rows with $\lambda_i$ boxes in the $i$-th row. The partition $\lambda$ is referred to as the shape of the diagram.
For instance, the Young diagram of shape $(4,3,2,2)$ is illustrated on the left in Figure \ref{4322}. In this paper, we introduce the notion of a semi-diagram, which is a Young
diagram with some shaded cells filled with a minus sign and subsequently turned into hollow cells. The diagram on the right in Figure \ref{4322} is a semi-diagram of shape $(4,3,2,2)$ with three minus signs.

To a semi-diagram $T$ of shape $\lambda$, we associate it with a weight as follows. For $1\leq i\leq k$, if the $i$-th row contains
$r_i$ shaded cells and $s_i$ minus signs, we define
its weight as $a_{r_i}z^{s_i}$. The weight
of a semi-diagram $T$, denoted by $w(T)$, is then defined to be the
product of weights of all rows. Of course,
the weight of an empty row
is set to be $a_0$. For example, the
weight of the semi-diagram in Figure \ref{4322} equals
$a_3a_2^2a_1z^3$.

In terms of semi-diagrams, the substitution \eqref{ai} can be  interpreted as filling some of the shaded cells of a Young
diagram  with a minus sign and turning them into hollow cells. Consider
only one row with $i$ shaded cells, which has
weight $a_i$. For $0\leq j \leq i$, there are ${i \choose j}$ ways to turn this
 row into a semi-diagram by placing  a minus sign to $j$ shaded cells and turning them
 into hollow cells. Any of the resulting semi-diagrams has weight $a_{i-j}z^j$.
 This operation is in accordance with the substitution \eqref{ai}.

For example, under the shear transformation, $a_3$ becomes \[
a_3'=a_3+3a_2z+3a_1z^2+a_0z^3,\]
 which
is the sum of the weights of the semi-diagrams
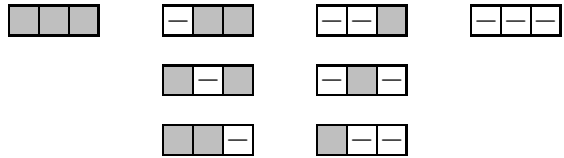
\[ \ytableaushort{\ \ \ }* [*(lightgray)]{3} \qquad \ytableaushort{- {*(lightgray)\ } {*(lightgray)\ } ,\none,{*(lightgray)\ }-{*(lightgray)\ }, \none,{*(lightgray)\ }{*(lightgray)\ }-} \qquad \ytableaushort{- - {*(lightgray)\ } ,\none,- {*(lightgray)\ }- , \none,{*(lightgray)\ }- -} \qquad \ytableaushort{- - -}
\]

We are now ready to prove Theorem \ref{SI}.

\noindent
{\it Proof of Theorem \ref{SI}.}
It suffices to show that for any monomial $a_\lambda$,
\begin{equation}\label{t0}
a'_\lambda=a'_{\lambda_1} a'_{\lambda_2} \cdots a'_{\lambda_k}
= \sum_{i=0} ^m D^i(a_\lambda)\frac{z^i}{i!}.
\end{equation}

 Let $K_i(\lambda)$ denote the set of semi-diagrams of
  $\lambda$ containing $i$ minus signs, where $0\leq i\leq m$.
 Set
 \begin{equation}
 \sum_{T\in K_i(\lambda)} w(T)=W_i(\lambda).
 \end{equation}
 Then we have
 \begin{equation}\label{t1}
 a'_\lambda= \sum_{i=0} ^m W_i(\lambda).
 \end{equation}
 In particular, the semi-diagrams containing only one minus sign give rise to the coefficient
 of $z$ in $I(a'_0,a'_1,\ldots,a'_n)$. This operation of filling only one shaded cell of a Young diagram with a
 minus sign and turning it into a hollow cell can be
 described by the action of the operator $D$  on $a_\lambda$, and this explains where the
 operator $D$ comes from combinatorially.

 Moreover, one realizes that in general the coefficient of $z^i$ in $I(a'_0,a'_1,\ldots,a'_n)$ can also be expressed in terms of the operator $D$.
 Indeed, the action of $D^i$ on $a_\lambda$ can be interpreted as placing
 $i$ distinguishable minus signs in a Young diagram and turning them
  into hollow cells. But for the coefficient of $z^i$ in $I(a'_0,a'_1,\ldots,a'_n)$, the minus signs in a
 semi-diagram are regarded indistinguishable. This yields the relation
 \begin{equation}\label{t2}
D^i(a_\lambda)\frac{z^i}{i!}=W_i(\lambda).
\end{equation}
Combining \eqref{t1} and \eqref{t2} gives \eqref{t0}. This completes the proof.
\qed

Now we turn to  the vertical shear transformation:
 \[ x= x'', \quad  y=zx''+y'',\]
 where $z$ is considered as a variable.
Under this transformation, the binary form $f(x,y)$ becomes
\begin{align*}
f''(x'',y'')&=a_0''{x''}^n+{n\choose1}a_1''{x''}^{n-1}y''+\cdots+a_n''{y''}^n,
\end{align*}
where, for  $0\leq i\leq n$,
\begin{equation}\label{verai}
a_i''=a_i+{n-i\choose1} a_{i+1}z+{n-i\choose2} a_{i+2}z^2+\cdots+a_n z^{n-i}.
\end{equation}

A polynomial $I(a_0,\,a_1,\ldots,\,a_n)$ with rational coefficients is called a semi-invariant with respect to the vertical shear transformation provided that
\begin{equation}
I(a_0,\,a_1,\ldots,\,a_n)=I(a''_0,\,a''_1,\ldots,\,a''_n).
\end{equation}
A polynomial $I(a_0,\,a_1,\ldots,\,a_n)$ is a semi-invariant with
respect to the vertical shear transformation if and only if  $\Delta (I)=0$, see Hilbert \cite{Hilbert93}.

To give a combinatorial interpretation of the operator $\Delta$, we need to
have a full picture of the Young diagram of a partition $\lambda$ contained in a $k \times n$ rectangle.
More precisely, we shall use shaded cells for the cells in the shape of $\lambda$ and
use hollow cells for the cells outside the shape of $\lambda$.
For example, below is the
depiction  of the partition $\lambda=(4,2,1,0)$ contained in a $4\times 5$ rectangle:
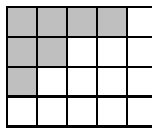
\[
\ytableaushort{\none,\none, \none,\none}
* {5,5,5,5}* [*(lightgray)]{4,2,1}
\]

Consider a single row diagram
 with $i$ shaded cells and $n-i$ hollow cells, whose weight is $a_i$.
For $0\leq j \leq n-i$, there are ${n-i \choose j}$ ways to turn this
 row into a semi-diagram by placing  a plus sign to $j$ hollow cells and turning them
 into shaded cells. Any of the resulting semi-diagrams has weight $a_{i+j}z^j$, where
 we define the weight of a plus sign to be $z$.
For example, for $n=6$, under the vertical shear transformation,
 $a_3$ becomes \[ a_3''=a_3+3a_4z+3a_5z^2+a_6z^3,\] which
is the sum of the weights of the semi-diagrams
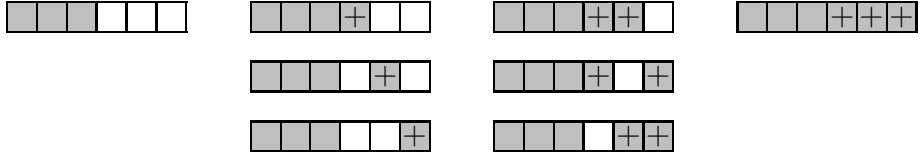
\vspace{7pt}
\[
\ytableaushort{\none}
*{6}* [*(lightgray)]{3}\qquad\ytableaushort{{*(lightgray)\ } {*(lightgray)\ } {*(lightgray)\ } {*(lightgray)+} \ \ ,\none,{*(lightgray)\ } {*(lightgray)\ } {*(lightgray)\ } \ {*(lightgray)+} \ ,\none,{*(lightgray)\ } {*(lightgray)\ } {*(lightgray)\ } \ \ {*(lightgray)+} }\qquad\ytableaushort{{*(lightgray)\ } {*(lightgray)\ } {*(lightgray)\ } {*(lightgray)+} {*(lightgray)+} \ ,\none,{*(lightgray)\ } {*(lightgray)\ } {*(lightgray)\ } {*(lightgray)+} \ {*(lightgray)+} ,\none,{*(lightgray)\ } {*(lightgray)\ } {*(lightgray)\ } \ {*(lightgray)+} {*(lightgray)+} }\qquad\ytableaushort{{*(lightgray)\ } {*(lightgray)\ } {*(lightgray)\ } {*(lightgray)+} {*(lightgray)+} {*(lightgray)+} }
\]
\vspace{1pt}

In particular, the action of the operator $\Delta$  can be interpreted
 in terms of the operation of filling only one  hollow cell in a Young diagram of  $\lambda$ with a plus sign and turning it into a shaded cell.
Analogous to Theorem \ref{SI}, we have the following expansion.

\begin{theo}\label{Rr}
For $n\geq 0$ and for any polynomial $I(a_0, a_1, \ldots, a_n)$ over the rational numbers, we have
\begin{equation}
I(a''_0, a''_1, \ldots, a''_n)=\sum_{i\geq0} \Delta^i I(a_0, a_1, \ldots, a_n)\frac{z^i}{i!}.
\end{equation}
\end{theo}

\section{Sylvester's Theorem}

The following   theorem of Sylvester \cite{Sylvester78b} establishes
a connection between the Gaussian coefficients
and semi-invariants.
For $0\leq m\leq nk/2$, let
\begin{equation}
\delta(k,n,m)=p(k,n, m)-p(k,n, m-1),
\end{equation}
with the convention that $p(k,n,-1)=0$.

\begin{theo} \label{weishu}
For $n,k\geq 0$ and $0\leq m\leq nk/2$, the number of semi-invariants of a binary $n$-form of degree $k$ and weight $m$ equals $\delta(k,n,m)$.
\end{theo}

 Theorem \ref{weishu} also takes the following form,  as Sylvester \cite{Sylvester78b}
 chose to work with.

\begin{theo}\label{nem}
For $n,k\geq 0$ and $0\leq m\leq nk/2$, the number of semi-invariants of a binary $n$-form of degree $k$ and weight not exceeding $m$ equals $p(k,n,m)$.
\end{theo}

The following identity \eqref{resI} of Cayley \cite{Cayley56} crops up in Sylvester's proof of the unimodality of the Gaussian coefficients, see  \cite{Sylvester78b}.
Cayley used this relation to construct covariants.

It might be worth mentioning that from a combinatorial point of view, the identity of Hilbert is easier
to justify than the identity of Cayley
for it does not involve the conditions on semi-invariants.
It should also be noted that to pass from Hilbert's identity
to Cayley's identity, the condition $m\leq nk/2$ is required
because it is a constraint for semi-invariants, see Cayley \cite{Cayley56}.

\begin{theo}\label{1i-1}
For $n,k\geq 0$, and $0\leq m\leq nk/2$, let $I$ be a semi-invariant of a binary $n$-form of degree $k$ and weight $m$, and let $c=nk-2m$.
Then, for $i\geq 1$,
\begin{equation}\label{resI}
{D}{\Delta}^i (I)
=i(c-i+1)\Delta^{i-1} (I).
\end{equation}
\end{theo}

Since $D(I)=0$ for a semi-invariant $I$, we see that the above relation is a consequence of the identity \eqref{kf} of Hilbert \cite{Hilbert93} with $a_\lambda$ being replaced by a
semi-invariant $I$.

To make the paper self-contained, we give an exposition of Sylvester's proof. The consideration of the dimension identity
\eqref{sum} makes the argument easier to understand in terms of an equality, instead of
arguing with two inequalities in opposite directions as described by Sylvester.

We define  $Q_n(k,m)$ as the
vector space of polynomials in $a_0, a_1, \ldots, a_n$
over the rational numbers that are homogeneous of degree $k$ and weight $m$.
We shall use $S_n(k,m)$
to denote the vector space of semi-invariants of degree $k$ and
weight $m$, that is,
\begin{equation} \label{Vkp}
S_n(k,m) = \{I \in Q_n(k,m)  \mid D (I) = 0 \}.
\end{equation}
The number of
semi-invariants of degree $k$ and weight $m$ of
a binary $n$-form is referred to as
the dimension of the vector space $S_n(k,m)$.
For example, $\dim S_4(4,6)=2$. Below are  two semi-invariants of degree 4
and weight 6 of a binary $4$-form:
\begin{align}\label{I1I2}
I_1&=3a_1^2 a_2^2-4a_1^3 a_3-2a_0 a_1 a_2 a_3
+3a_0^2 a_3^2+4a_0a_1^2a_4-4a_0^2a_2a_4,\\[6pt]
I_2&=a_0 a_2^3-2a_0 a_1 a_2 a_3+a_0^2 a_3^2
+a_0 a_1^2 a_4-a_0^2 a_2 a_4.\label{I1I2-b}
\end{align}
Notice that for $m=0$, $\dim S_n(k,0)=p(k,n,0)=1$.

\noindent
{\it Proof of Theorem \ref{nem}.}
For $0\leq i\leq m+1$, let \[ V_i=D^i(Q_n(k,m)), \]
and let
\[ T_i\colon V_{i-1}\rightarrow V_i,\]
where $1\leq i\leq m+1$, that is,
$T_i (I) = D(I)$ for any $I \in  V_{i-1}$.
Then, the kernel of $T_i$, denoted by $\ker T_i$, namely,
\[ \ker T_i= \{ I\in V_{i-1} \mid D(I)=0\} , \]
is a subspace of $V_{i-1}$. Hence
\begin{equation}\label{sum}
\dim V_{i-1}=\dim\ker T_i+\dim V_i.
\end{equation}
Notice that while acting on a monomial in $a_0, a_1, \ldots, a_n$,
the operator $D$ preserves the degree and lowers the weight by one.
On the other hand,
 $Q_n(k,0)$ is generated by $a_0^k$, which is a semi-invariant, and so
 $D(Q_n(k,0))=0$. It follows that \[ V_{m+1}=D^{m+1}(Q_{n}(k,m))=0.\]
Iterating \eqref{sum} gives
\begin{equation}\label{eq}
\dim V_0=  \dim \ker T_1 + \dim \ker T_2 + \cdots +  \dim \ker T_{m+1}.
\end{equation}
It is apparent that
\begin{equation}\label{sub}
\ker T_i\subseteq S_n(k,m-i+1).
\end{equation}
The real challenge is to show that
\begin{equation}\label{sub=}
\ker T_i =  S_n(k,m-i+1),
\end{equation}
that is, the successive applications of the operator $D$ resulting in $V_{i-1}$ do not leave out any semi-invariants in $S_n(k,m-i+1)$.

For $i=1$, nothing needs to be said since by definition,
\begin{equation}
\ker T_1=S_n(k,m).
\end{equation}
But for $i=2$,
what does \eqref{sub=} mean? Note that a semi-invariant $I$
in $S_n(k,m-1)$ should come from a polynomial in $Q_n(k,m-1)$. However, \eqref{sub=}
 indicates
that we can restrict our attention only to $D(Q_n(k,m))$, which is a subspace
of $Q_n(k,m-1)$, and we can still get all the semi-invariants in $S_n(k,m-1)$.

Sylvester realized that
in some sense the operator $D$ is the inverse of the
operator $\Delta$, as guaranteed by the identity \eqref{resI} of Cayley. In other words, the action of the operator $\Delta$ ensures that every
semi-invariant in $S_n(k,m-i+1)$ can be shielded from the annihilation of the operator $D$.

For example, for $n=3, k=3$ and $m=0$,
it is clear that $Q_3(3,0)$ is generated by $a_0^3$, which is a semi-invariant of degree three and weight zero. Now, $V_3 = D^3(Q_{3}(3,3))$.
One may wonder whether  $a_0^3$
 is still there in $V_3$. Employing the operator $\Delta$, we find that
\[\Delta^3 (a_0^3)=18a_0^2a_3+324a_0a_1a_2+162a_1^3,\]
which is a polynomial in $Q_3(3,3)$.
Then it is easily verified that
\[D^3 \Delta^3 (a_0^3) =3024a_0^3. \]
So $a_0^3$ remains in $V_3$, and this is in accordance with the fact that $S_3(3,0)$ is generated by $a_0^3$.

The above reasoning is valid for the general case. For any semi-invariant $I$ in $S_n(k,m-i+1)$,
$\Delta^{i-1}(I)$ falls into $Q_n(k,m)$. Thanks to the identity \eqref{resI}, we
deduce that by successively applying the operator $D$  to  $\Delta^{i-1}(I)$,  one recovers the semi-invariant $I$ if we do not mind the nonzero constant. To be more specific, we deduce that $I$ truly belongs to $V_{i-1}$,
and hence the proof is complete.
\qed

Examining the proof of Sylvester, one sees that what Sylvester tried to
 demonstrate is the following property of $D$.

\begin{theo}
For $n,k\geq 0$ and $1\leq m\leq nk/2$, we have
\begin{equation}\label{sur}
Q_n(k,m-1)= D(Q_n(k,m)),
\end{equation}
or equivalently,
the transformation $D$ is a surjection from
$Q_n(k,m)$ to $Q_n(k,m-1)$.
\end{theo}

Once the above surjectivity is in hand, it immediately follows that
the number of semi-invariants, namely, the dimension of the kernel of $D$,
is given by $\delta(k,n, m)$.
As far as the unimodality of the Gaussian coefficients is concerned,
we see that Sylvester's proof  also contains a justification of the injectivity of the transformation $\Delta$.

\begin{theo}
For $n,k\geq 0$ and $1 \leq m\leq nk/2$, the transformation $\Delta$ is an
injection from $Q_n(k,m-1)$ to $Q_n(k,m)$.
\end{theo}

Proctor \cite{Proctor82} came up with a proof of the unimodality of the Gaussian coefficients by
introducing two operators different from $D$ and $\Delta$, while taking
a notice of the surjectivity and injectivity of $D$ and $\Delta$. It would be appealing to
reach a better understanding of these properties from a combinatorial
angle.

\section{A Combinatorial Proof of Hilbert's Identity}

Based on the combinatorial interpretations of the  operators $D$ and $\Delta$, we give  a combinatorial proof of the identity  of Hilbert  \cite{Hilbert93},  as stated below.

\begin{theo}
\label{keyf}
For $n,k\geq 0$, and $0\leq m\leq nk$, let $\lambda$ be a partition of $m$ contained in a $k\times n$ rectangle, and let $c=nk-2m$. Then, for $i\geq1$,
\begin{equation}\label{kf-s-2}
D\Delta^i(a_\lambda)-\Delta^i D(a_\lambda)=i(c-i+1)\Delta^{i-1}(a_\lambda).
\end{equation}
\end{theo}

The above identity of Hilbert plays a fundamental role in the characterization of invariants as well as the construction of covariants.

To present a combinatorial interpretation of the above relation, we recall
that the diagram of a partition $\lambda$ contained in a $k\times n$ rectangle
 contains  shaded cells inside the shape of $\lambda$ and hollow cells
  outside the shape of $\lambda$. A semi-diagram will also be represented in
  the same manner. While we shall encounter some signs filled
  in semi-diagrams, the weight of a row with $r$ shaded cells will be defined by
  $a_r$. So the signs will not affect the weight of a semi-diagram.

\noindent{\it Proof of Theorem \ref{keyf}.}
The action of $D\Delta^i$ on $a_\lambda$ can be understood as placing $i$ distinguishable plus signs in hollow cells of the Young diagram of  $\lambda$ and turning them into
 shaded cells, and then placing a minus sign in a shaded cell and turning it into a hollow cell.

Similarly, the application of $\Delta^i D$ to $a_\lambda$ means to fill a shaded
 cell of the Young diagram of shape $\lambda$ with a minus sign and to turn it into
  a hollow cell,  then fill $i$ hollow cells with distinguishable plus signs and turn them
  into shaded cells.

Now, we use the symbol $\pm$ to denote a cell that is filled with a plus sign first and subsequently filled with a minus sign. Similarly, $\mp$ denotes a cell that is filled with a minus
sign first and subsequently filled with a plus sign.

With respect to the computation of $D \Delta^i (a_\lambda) - \Delta^i D (a_\lambda)$,
the semi-diagrams that do not contain any $\pm$ or $\mp$ signs would cancel out.

In fact, it is readily seen that there is a one-to-one correspondence between the semi-diagrams in $D \Delta^i (a_\lambda) $ that do not contain any $\pm$ signs  and the semi-diagrams in $\Delta^i D (a_\lambda)$ that do not contain any $\mp$ signs.
For example, when $k=4, n=5$ and $\lambda=(4,2,1,0)$, the following semi-diagram occurs in $D \Delta^3 (a_\lambda)$, and it is also in $\Delta^3 D (a_\lambda)$:
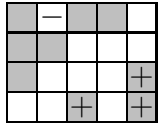
\[
\ytableaushort{\ {*(white)-}\ \ \ ,\none, \ \ \ \ {*(lightgray)+},\ \ {*(lightgray)+}\ {*(lightgray)+}}* {5,5,5,5}
* [*(lightgray)]{4,2,1}
\]

Therefore, it suffices to consider the cases when the sign $\pm$ or $\mp$ is
  involved.

Let us consider the semi-diagrams generated by $D\Delta^i (a_{\lambda})$
containing the sign $\pm$.
For example, when $k=4, n=5$ and $\lambda=(4,2,1,0)$,
the following semi-diagram is an illustration of such a configuration generated by
$D \Delta^3 (a_\lambda)$:
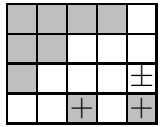
\[
\ytableaushort{ \none,\none, \ \ \ \ \pm,\ \ {*(lightgray)+}\ {*(lightgray)+}}* {5,5,5,5}
* [*(lightgray)]{4,2,1}
\]
These configurations under consideration can  be produced from the
Young diagram of $\lambda$ in a $k \times n$ rectangle by
placing $i$ distinguishable plus signs in  hollow cells  and subsequently adding a minus sign to a shaded cell with a plus sign and turning
it into a hollow cell. On the other hand, semi-diagrams with a $\pm$ cell and $i-1$ distinguishable plus signs  (outside
 the shape of $\lambda$, to be precise) can  be constructed in an alternative way.

We may choose $i-1$ plus signs from the $i$ distinguishable plus signs,
  place them in the hollow cells of the Young diagram of  $\lambda$ contained
 in a $k\times n$ rectangle, turn them into shaded cells, and finish with filling
  a hollow cell with the $\pm$ sign and keeping it hollow. Notice that the location of a
 hollow cell with the $\pm$ sign does not affect the weight of a semi-diagram.
 That is to say, the cell with the symbol $\pm$ might as well be viewed just as a hollow cell.
 Therefore, if the $\pm$ sign is not taken into consideration,
  the semi-diagrams containing $i-1$ distinguishable plus signs are generated by
applying the operator $\Delta^{i-1}$ to $a_\lambda$.
Note that there are $nk-m-(i-1)$  hollow cells left for the moment, any of which
can be chosen as a residence of the $\pm$ sign.
It follows that the total contribution of weights in this case amounts to
\begin{equation}
\label{pm-1}
i(nk-m-(i-1))\Delta^{i-1}(a_{\lambda}).
\end{equation}

We now consider the semi-diagrams generated by $\Delta^i D(a_{\lambda})$
containing the sign $\mp$. For example,
when $k=4, n=5$ and $\lambda=(4,2,1,0)$,
the following semi-diagram  arises in $\Delta^3 D (a_\lambda)$:
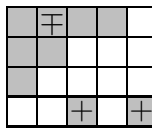
\[
\ytableaushort{\ \mp \none,\none,\none,\ \ {*(lightgray)+}\ {*(lightgray)+}}* {5,5,5,5}
* [*(lightgray)]{4,2,1}
\]
In general, the involved  configurations in this case
 can be generated from a Young diagram of shape $\lambda$ contained
in a $k\times n$ rectangle by placing
a minus sign in a shaded cell, and then distributing $i$ distinguishable plus
 signs to $i-1$ purely hollow cells  and the remaining  plus sign to the cell with a minus sign. Note that the cells carrying the $\mp$ sign as well as the plus sign are all shaded
 in the end.

There is also another way to construct such configurations.
We may choose $i-1$ plus signs from the $i$ distinguishable plus signs,  place
 them in the hollow cells and turn them into shaded cells, and then fill a shaded cell with the sign $\mp$ and keep it shaded. Since there are $m$ choices
  for a shaded cell to be filled with the sign $\mp$, the
  weights of all the feasible configurations generated by
   $\Delta^i D (a_{\lambda})$ add up to
\begin{equation}
\label{pm-2}
i m\Delta^{i-1}(a_{\lambda}).
\end{equation}

Casting up \eqref{pm-1}  and \eqref{pm-2}, we arrive at
\begin{equation}
i(nk-2m-i+1)\Delta^{i-1}(a_{\lambda}),
\end{equation}
in agreement with the right hand side of \eqref{kf-s-2}. This completes the proof.
\qed

It is worth mentioning that Hilbert obtained another identity on $D$ and $\Delta$. As before, let $c=nk-2m$. Then, for $i\geq 1$,
\begin{equation}\label{i1}
D^i\Delta(a_\lambda)-\Delta D^i (a_\lambda)=i(c+i-1)D^{i-1}(a_\lambda).
\end{equation}
Using a similar approach to the identity \eqref{kf-s-2}, it is not hard to
provide a combinatorial
proof of \eqref{i1}. The detailed description is left out.
 Utilizing  \eqref{i1},  Hilbert demonstrated that the number of invariants of a binary $n$-form of degree $k$ and weight $m=nk/2$ equals $\delta(k,n,nk/2)$, where at least one of $k$ and $n$ is even.
In fact, Hilbert deduced that the operator $\Delta$ is an injection from $Q_n(k,nk/2-1)$ to $Q_n(k,nk/2)$.
If not, there would exist a nonzero polynomial $I$ that can be written as a linear combination of the basis elements of $Q_n(k,nk/2-1)$ such that in some way the action of $\Delta$ on $I$ yields an invariant in $Q_n(k,nk/2)$, that is,  $D\Delta(I)=0$.
On the other hand, there must exist $i$ such that $ D^i(I)=0$ but $ D^{i-1}(I)\neq0$.
 This incurs a contradiction with \eqref{i1}.

We also remark that the argument of Hilbert is valid for semi-invariants with $0 \leq m \leq nk/2$ where the parity constraint on $n$ and $k$ may be lifted. So it can be regarded as an alternative proof of Sylvester's theorem.
Note that in the case when at least one of $k$ and $n$ is even and $m=nk/2$,
a semi-invariant of degree $k$ and weight $m$ turns out to be an invariant.

\section{The Additivity Lemma of Pak and Panova}

In this section, we present a derivation of the additivity lemma
of Pak and Panova in the context of semi-invariants. As will be seen,
it is a consequence of the ring property of semi-invariants.
In view of Theorem \ref{weishu}, to prove the strict unimodality of $n+k_1+k_2\brack n$, it suffices to show that for any $2\leq m\leq \lfloor n(k_1+k_2)/2\rfloor$, there exists a semi-invariant of degree $k_1+k_2$ and weight $m$.

\noindent
{\it Proof of Lemma \ref{add}.}
It is not difficult to see that
under the conditions on $n, k_1, k_2$, for
any $4\leq m\leq \lfloor n(k_1+k_2)/2\rfloor$, we can always express $m$ as $m_1+m_2$ such that $2\leq m_1\leq \lfloor nk_1/2\rfloor$ and $2\leq m_2\leq \lfloor nk_2/2\rfloor$.

Since  at least one of $k_1$, $k_2$ and $n$ is even, we find that
\begin{equation}\label{haveeven}
\lfloor nk_1/2\rfloor+\lfloor nk_2/2\rfloor=\lfloor n(k_1+k_2)/2\rfloor.
\end{equation}
This takes care of the case $m=\lfloor n(k_1+k_2)/2\rfloor$.

For any $\lfloor nk_1/2\rfloor+2\leq m < \lfloor n(k_1+k_2)/2\rfloor$, taking
$m_1=\lfloor n k_1/2\rfloor$, a simple computation shows that the corresponding
$m_2$ falls into the right range, that is, $2\leq m_2< \lfloor nk_2/2\rfloor$.

For any $4\leq m < \lfloor nk_1/2\rfloor+2$, we may choose $m_2$ to be $2$ and set $m_1=m-2$. It can be checked that  $2\leq m_1< \lfloor nk_1/2\rfloor$.

We now assume that $4\leq m \leq \lfloor nk_1/2\rfloor+\lfloor nk_2/2\rfloor$
and $m=m_1+m_2$, where  $2\leq m_1\leq \lfloor nk_1/2\rfloor$ and $2\leq m_2\leq \lfloor nk_2/2\rfloor$. By the assumptions for  ${n+k_1\brack n}$ and ${n+k_2\brack n}$, we see that there exists a semi-invariant $I$ of a binary $n$-form of degree $k_1$ and weight $m_1$  and a semi-invariant $J$ of a binary $n$-form of degree $k_2$ and weight $m_2$.
Consequently, $IJ$ is a semi-invariant of a binary $n$-form of degree $k_1+k_2$ and weight $m$.

We now turn to the remaining cases $m=2$ and $m=3$.
Since $k_1,k_2,n\geq 2$ and at least one of $k_1$, $k_2$ and $n$ is greater than two,
it is evident that at least one of $nk_1/2$ and $nk_2/2$ is greater than or equal to three.
Let us assume that $nk_1/2\geq 3$.
By the assumption for ${n+k_1\brack n}$, there is a semi-invariant of degree $k_1$ and weight $m_1$ for any $2\leq m_1\leq nk_1/2$.
For $m_1=2$, suppose that $I$ is a semi-invariant of degree $k_1$ and weight two.
For $m_1=3$, suppose that $J$ is a semi-invariant of degree $k_1$ and weight three.
Clearly, $a_0^{k_2}$ is a semi-invariant of
degree $k_2$ and weight zero.
It follows that $a_0^{k_2}I$ is a semi-invariant of degree $k_1+k_2$ and weight two and $a_0^{k_2}J$ is a semi-invariant of degree $k_1+k_2$ and weight three.
This completes the proof.
\qed

\vspace{0.5cm}
\noindent{\bf\large Acknowledgments.}\ \
We wish to thank the referee for invaluable comments.
This work was done under the auspices of the National Science Foundation of China.

\end{document}